\documentclass{amsart}
\usepackage{amsfonts,amsmath,amsthm,amssymb}
\usepackage[latin1]{inputenc}
\usepackage{graphicx}
\usepackage{url}

\usepackage{fancyhdr}

\usepackage[breaklinks]{hyperref}
\usepackage[mathscr]{euscript}
\usepackage{enumerate,multicol,array}

\makeatletter
\def\newrefformat#1#2{%
  \@namedef{pr@#1}##1{#2}}
\def\prettyref#1{\@prettyref#1:}
\def\@prettyref#1:#2:{%
  \expandafter\ifx\csname pr@#1\endcsname\relax%
    \PackageWarning{prettyref}{Reference format #1\space undefined}%
    \ref{#1:#2}%
  \else%
    \csname pr@#1\endcsname{#1:#2}%
  \fi%
}
\makeatother

\newrefformat{eq}{\textup{(\ref{#1})}}
\newrefformat{cha}{Chapter \ref{#1}}
\newrefformat{sec}{Section \ref{#1}}
\newrefformat{tab}{Table \ref{#1} on page \pageref{#1}}
\newrefformat{fig}{Figure \ref{#1} on page \pageref{#1}}

\makeatletter
\def\indsym#1#2{%
  \setbox0=\hbox{$\m@th#1x$}%
  \kern\wd0%
  \hbox to 0pt{\hss$\m@th#1\mid$\hbox to 0pt{$\m@th#1^{#2}$}\hss}%
  \lower.9\ht0\hbox to 0pt{\hss$\m@th#1\smile$\hss}%
  \kern\wd0} 
\def\nindsym#1#2{%
  \setbox0=\hbox{$\m@th#1x$}%
  \kern\wd0%
  \hbox to 0pt{\mathchardef\nn="3236\hss$\m@th#1\nn$\kern1.4\wd0\hss}
  \hbox to 0pt{\hss$\m@th#1\mid$\hbox to 0pt{$\m@th#1^{#2}$}\hss}%
  \lower.9\ht0\hbox to 0pt{\hss$\m@th#1\smile$\hss}%
  \kern\wd0}

\makeatother

 \def\bw{\bar w}
\def\p{\pi}

\def\ben{\begin{enumerate}}\def\een{\end{enumerate}}
\def\bdc{\begin{description}}\def\edc{\end{description}}
\def\bitm{\begin{itemize}}\def\eitm{\end{itemize}}
\def\bdf{\begin{defin}}\def\edf{\end{defin}}
\def\bth{\begin{theo}}\def\eth{\end{theo}}
\def\bfc{\begin{fact}}\def\efc{\end{fact}}
\def\bco{\begin{coro}}\def\eco{\end{coro}}
\def\brm{\begin{rem}}\def\erm{\end{rem}}
\def\blm{\begin{lemma}}\def\elm{\end{lemma}}
\def\bnt{\begin{nota}}\def\ent{\end{nota}}
\def\bex{\begin{exe}}\def\eex{\end{exe}}

\def\txtt{\textit}

     \def\Sum{\sum}

\def\O+{\oplus}

\newtheorem{theo}{Theorem}[section]

\newtheorem{coro}[theo]{Corolary}
\newtheorem{lemma}[theo]{Lemma}

\theoremstyle{definition}
\newtheorem{defin}[theo]{Definition}
\newtheorem{fact}[theo]{Fact}
\theoremstyle{remark}
\newtheorem{exe}[theo]{Example}
\newtheorem{rem}[theo]{Remark}
\newtheorem{nota}[theo]{Notation}

\newcommand{\R}{\ensuremath{\mathbb{R}}}

\author{Francisco Gutierrez}
\address{Francisco Gutierrez San\'in
\\ Universidad Sergo Arboleda
\\ Departamento de Matem\'aticas}

\author{Camilo Argoty}
\address{Camilo Argoty
\\ Universidad Sergo Arboleda
\\ Departamento de Matem\'aticas
\\ }

\author{Stefany Moreno}
\address{Stefay Moreno
\\ Universidad Sergo Arboleda
\\ Departamento de Matem\'aticas
\\ }

\title{A generalization of a classical model in contract theory: The agent behavior}
\date{}
\begin{document}

\maketitle

\begin{abstract}
We present a first approximation of agent behaviour in a generalized model in contract theory. This model relaxes some of the the assumptions of one of the classical models allowing to include a broader range of agents. We introduce the motivation for the agent and reinterpret the classical definition of risk perception. Besides, we analyze different scenarios for the relation between the effort exerted by the agent and the probability that he gets an especfic result.  
\end{abstract}

\section{Introduction}

\subsection{The classical model}
A contract is an agreement between two or more parties that defines a set of mutual obligations for them. In the classical model the contractual relation is established between two parties: the principal and the agent \cite{Bolton}. Each party has to make one decision: The principal decides the wage $w$ that the agent will receive according to the monetary result $x$ obtained from the contractual relation, and the agent decides how much effort $e$ he will exert. It is worthy to note that the monetary result $x$ does not depend only on the effort exerted by the agent but also on a random variable associated with all possible external conditions that affect the final result and are not controlled by any of the parties \cite{Bolton}. 

The relation between the variables $e$, $x$ and $w$ and the utility of each party is defined by an {\it utility function}. This function is specific for each player and expresses his preferences with respect to the risk. The utility function for the principal is $B(x-w)$. It is assumed that $B'>0$, which means that the utility of the principal increases with $x$ and decreases with $w$ and $B''\leq0$ which means that it is concave \cite{Macho}. The second derivative will determine the degree of risk aversion of the principal: if $B''<0$ he is risk averse, if $B''=0$ he is risk neutral and if $B''>0$ he is risk-seeking.

On the other hand, the utility function for the agent is defined by $U(w,e)=u(w)-v(e)$. It is assumed that, $u'(w)>0$, $u''(w)\leq0$, $v'(e)>0$ and $v''(e)\geq0$. The assumptions about $u(w)$ imply that the utility $u$ increases with $w$ and that the agent is either risk-averse or risk-neutral in terms of his payoff \cite{Macho}.  From $v(e)$, it means that $u$ decreases with $e$ with the marginal disutility of the effort not decreasing.

Using the Kuhn-Tucker Theorem to maximize the utility of the principal restricted to the agent's participation condition -which is having an expected utility higher than his reservation utility- it is possible to know how the contract should be designed according to the risk preferences of both parties. If one of the parties is risk averse the other should assume all the risk in the contractual relation, while if both are risk-averse they have to share the risk according to their degree of aversion \cite{Macho}.

\subsection{Limitations of the classical model}
First, note that it is assumed that the utility of the agent always increases with the wage received and always decreases with the effort exerted. This assumption is valid in the context of agents who make a repetitive work (e.g. machine operators), but it does not hold for all type of agents. It is possible that for an agent it is a loss not to exert any effort, i.e. the agent has an inner need of working. Even more, there may be an interval of efforts for which the agent has some utility in exerting them (e.g. volunteers, suicide bombers). Indeed, the initial motivation of this work was to extend the classical model of contract theory to the case of this kind of agents. 

Second, note that the other main assumption of the classical model is that the optimization is made only by one of the parties. Usually, that party is the principal who is the one designing the contract. As in the case of the Kuhn-Tucker maximization mentioned before, the agent just imposes some restrictions and the principal maximizes with respect to them. We consider that there is no reason to assume that just one party is maximizing. By doing the maximization for the agent first, it is possible to use other type of restrictions while doing the maximization for the principal. So, we replaced the participation and incentive compatibility restrictions for maximization restrictions for the agent.   

\section{Proposed model}
\subsection{The Principal-agent game}
We state the problem in the following way. We have two players: Agent $(A)$ and Principal $(B)$. Principal wants the agent to produce certain good $x$ among a set $X$ of possible goods $x_1,\dots, x_n$
\bdf
A \txtt{contract} is a function $w:X\to \R$. Equivalently, a contract is a vector $\bw=(w_1,\dots,w_n)$, such that $w_i=w(x_i)$
\edf

In the first move, principal chooses a contract $\bw$. Then, agent chooses a level of effort $e$ to exert at working. In the third step nature chooses a good $x_{i_0}$ according to a probability distribution $p_i(e)$ on $X$ that depends on the agent $A$ and the effort $e$ exerted by him, where $e \in [e_{min},e_{max}]$

It is worthy to analyze this probability distribution $p_i(e)$. This distribution is given by $p_i(e)=\{p_1(e),\dots,p_n(e)\}$, where $p_i(e)$ is the probability for the agent $A$ to produce the good $x_i$ when $A$ exerts an effort $e$. 

\bdf
The function $P(A)$ such that $P(A)=p_i(e)$ is called the \txtt{profile} of agent $A$
\edf

In the last step, principal and agent receive payments $\p_E(x,w)$ and $\p_A(w,e)$ respectively.

The goal of this section is solve this game by backward induction. We suppose by a first approach, that the distribution function $p_i(e)$ is known to both agent and principal.

In order to start this backward induction we need to analyze payment functions for agent and principal.

\subsection{Payment for agent}
We model payment function for agent by introducing two new functions $u(w)$ and $v(e)$ this way:
\[\p_A^{\bw}(w,e)=u(w)-v(e),\]
where, following tradition, $u(w)$ is an \txtt{utility} function depending on the wage $w$ received by the agent according to the contract $\bw$ selected by the principal, and $v(e)$ is a loss to the agent cause by the exertion of an effort $e$.

We propose some changes in function $v(e)$. First of all, we consider that it is not necessary to consider it a loss all the times. We assume that it is possible that for an agent it is a loss not to exert any effort, i.e. the agent has an inner need of working. Even more, there may be an interval of efforts for which the agent has some utility in exerting them. The usual agent become a particular case of this kind of function.
The main (perhaps the only) characteristic of $v$ is that $v^{\prime\prime}>0$. 


\bdf
Let $A$ be an agent with payment function $\p_A^{\bw}(w,e)=u(w)-v(e)$ and profile $P(A)$. Then the term $\Sum_{i=1,\dots,n}p_i(e)u(w_i)-v(e)$ is called the \txtt{payment expectation} for the agent $A$ and is denoted by $E_{A}^{\bw}(e)$. Note that the function $E_{A}^{\bw}$ depends on the contract selected by the principal ${\bw}$ and it has only one parameter: $e$. Once the principal selects the contract, the payment expectation function of the agent is established. Then, by modifying his effort -which is the only variable he can control-, the agent maximizes his expectation utility. 
 
Let $A$ be an agent, let $E_{A}^{\bw}(e)$ his benefit expectation function, and $\bw$ a contract. The derivative,
\[\frac{\partial E_A^{\bw}}{\partial e}\]
is called the \txtt{motivation function} of agent $A$ under contract $\bw$ and is denoted by $Mt_A^{\bw}(e)$.
Also, the derivative,
\[\frac{\partial Mt_A^{\bw}}{\partial e}\]
is called the \txtt{persistence function} of agent $A$ under contract $\bw$ and is denoted by $Prst_A^{\bw}(e)$. In case this function is negative, we can replace it by its opposite
\[-\frac{\partial Mt_A^{\bw}}{\partial e}\]
which is called the \txtt{transience function} of agent $A$ under contract $\bw$ and is denoted by $Tr_A^{\bw}(e)$
\edf

\subsection{The agent's problem}
Following backward induction, supposing that Principal $E$ has chosen a contract $\bw$, the agent has to solve the following optimization problem:
\[\max_{e\in [e_{min},e_{max}]}\bigl(\Sum_{i=1,\dots,n}p_i(e)u(w_i)-v(e)\bigr)\]
 In other words, 
 \[\max_{e\in [e_{min},e_{max}]} E_{A}^{\bw}(e)\]
 
 This maximization problem take us to three possible scenarios for $e^*$:
 
 1) $e^* \in (e_{min},e_{max})$
 
 By elementary calculus, the agent either chooses $e^*$ such that 
\[\frac{\partial E_A^{\bw}}{\partial e}=0\]
 which is equivalent to saying that the agent either chooses rejects the contract or chooses an $e^*$ such that $Mt_{A}^{\bw}(e^*)=0$ and
\[Prst_A^{\bw}(e^*)<0\text{ or equivalently,
}Tr_A^{\bw}(e^*)>0\]

2) $e^*=e_{min}$ which is  $E_{A}^{\bw}(e_{min}) \geq E_{A}^{\bw}(e) \text{ for all } e \in [e_{min},e_{max}]$ 

3)  $e^*=e_{max}$ which is $E_{A}^{\bw}(e_{max}) \geq E_{A}^{\bw}(e) \text{ for all } e \in [e_{min},e_{max}]$ 

Note that these scenarios are not exclusive i.e. it is possible that there is more than one maximum. 

\subsubsection{The Risk interpretation} We introduce one main shift to the usual way of treating risk: The position with respect to risk depends on the  utility expectation function of the agent, $E_{A}^{\bw}(e)$ and not only on his utility with respect to the payment $u(w)$. As it was explained before  $E_{A}^{\bw}(e)$ is defined by the contract. So, if the position with respect to risk depends on $E_{A}^{\bw}(e)$, this means that it is not an inherent characteristic of the agent but it changes depending on the type of contract offered by the principal. In other words, the same agent can have different positions with respect to risk for different type of contracts.\smallskip
 
 {\it Classical risk interpretation}
 
As it was mentioned before, in the classical risk interpretation the position of the agent with respect to risk depends on the sign of $u''$, which is constant for the agent. Let us consider each case, and its relation to the three maximization scenarios proposed above:
 
 (1) If the agent is risk- averse $Prst_A^{\tilde{w}}(e)<0$ for every value of $e$, so he is necessary in the first maximization scenario. This corresponds to classical agents that exert an amount of effort that is in $(e_{min},e_{max})$. In other words, they maximize their utility working (which means $e > e_{min}$) but not to his maximum potential ( $e > e_{max}$).
 
 (2) If the agent is risk-seeking $Prst_A^{\tilde{w}}(e)>0$ for every value of $e$, so it is not possible for him to find a maximum in $(e_{min},e_{max})$ which means that he is either in the second or the third case of the maximization problem. This implies that if the agent is always risk-seeking he will go for extreme values for the effort i.e. $e_{min}$ or $e_{max}$
 
 (3) If the agent is risk-neutral $Prst_A^{\tilde{w}}(e)=0$ which means that his expectation utility function is a straight line, which is classical result in contract theory. Note that the line can be an horizontal line, or it can have a slope different from 0. In the first case, the agent is indifferent to the al possible outcomes obtained by exerting any $e \in [e_{min},e_{max}]$. In the second case, it means that the only maximization scenarios possible for him are the second and the third. In other words, he can only find a maximum for his utility in  $\{e_{min}, e_{max}\}$. 

\section{Initial Results}

\subsection{Invisible effort}

The first case that we analyze is when is not possible for the principal to determine the effort exerted by the agent. This is a typical information asymmetry in contract theory and can be modeled by considering the function $p_i (e)$. Let us consider the more extreme case i.e. when the result obtained by the agent does not depend on his effort. In terms of  $p_i (e)$ this is $p_i (e)=p_i$ for every $e \in [e_{min},e_{max}]$. 

In this scenario, the utility expectation function is $E_{A}^{\bw}= \Sum_{i=1,\dots,n}p_iu(w_i)-v(e)$. As the first term is constant, differentiating by $e$ we get  $ Mt_A^{\bw}(e)= -v'(e)$, which means that in this scenario the agent will only minimize his disutility with respect to effort. Since in terms of the result the efforts are equivalent, the agent will choose the effort that implies the smallest disutility. So, he will be motivated to increase his effort just by an intrinsic need of working. Besides, we have  $Prst_A^{\bw}(e)= -v''(e)$. As it was explained before $v(e)$ is assumed to be a concave function. Then, $E_{A}^{\bw}=-v(e)$ is convex and  $ -v''(e)<0$. This implies that when the effort is invisible for the principal, the agent is risk-averse with respect to it. In other words, the agent will not take any risk by increasing the amount of effort exerted as he will not increase the probability of getting a better result. 

Note that the previous analysis can be extended to a case in which the function $p_i (e)$ is not a constant but it is almost independent of the effort (or changes very slightly with it). Actually, these are the type of functions that could exist in a real context.
 
\subsection{Two possible outcomes}

Let us consider the case in which there are just two possible results $x_1$ (good result) and $x_2$ (bad result) for the work exerted by the agent. As $p_1(e)+p_2(e)=1$ for any $e$, we can consider only the function $p_1(e)$. 

Let us assume that there is a linear relation between the probability of getting $x_1$ and the effort exerted by the agent. Then, $p_1 (e)=Ce +h$ with $C$ and $h$ such that $0\leq p_1(e_{min}), p_1(e_{max})\leq1$. Then, the utility expectation function is $E_{A}^{\bw}= \Sum_{i=1}^2 p_i(e)u(w_i)-v(e)$. Differentiating by $e$ we get $Mt_A^{\tilde{w}}(e)= C(u(w_1)-u(w_2))-v'(e)$. This implies that the there is a maximum when $ C(u(w_1)-u(w_2))=v'(e)$. Since it is assumed that the function $v'(e)$ is concave, its derivative increases when $e$ increases. From the previous equation, $C$ and $u(w_1)-u(w_2)$ are proportional to $v'(e)$ and then hold a direct relation with $e$. If $C$ increases, the probability of getting the good result by increasing the effort is higher. This explains why there is a motivation for the agent to exert a higher effort. On the other side, if there is a ig difference between the utility received with each one of the wages, it is also profitable for him to increase his effort.  Besides, $Prst_A^{\tilde{w}}(e)= -v''(e)$ Since, $-v''(e)<0$ the agent is risk-averse with respect to $e$.

Note that this means that in terms of risk aversion this agent is equivalent to the one in a situation of invisible effort. His risk perception only involves the disutility for exerting $e$ and not the potential utility generated by $w$. Then, in these scenarios the agent will always be risk-averse. Indeed, the only way in which an agent can be risk-seeking is if there is a non linear relation between $e$ and $p$.

\section{Conclusions}
It was shown that, by generalizing the classic model of contract theory, it is possible to include a broader type of agents and still get the same conclusions about the effect of risk preferences in the contractual relation. Moreover, we showed that the agent's motivation is a key factor in the contractual relation that is usually ignored in the classical contract theory model. The future direction of this work consists in including information asymmetries between the principal and the agent. This will be done by using information theory to deal with the unknown information in the contractual relation.

\end{document}